# Minimal Input Selection for Robust Control

Zhipeng Liu, Yao Long, Andrew Clark, Phillip Lee, Linda Bushnell, Daniel Kirschen, and Radha Poovendran

*Abstract*— This paper studies the problem of selecting a minimum-size set of input nodes to guarantee stability of a networked system in the presence of uncertainties and time delays. Current approaches to input selection in networked dynamical systems focus on nominal systems with known parameter values in the absence of delays. We derive sufficient conditions for existence of a stabilizing controller for an uncertain system that are based on a subset of system modes lying within the controllability subspace induced by the set of inputs. We then formulate the minimum input selection problem and prove that it is equivalent to a discrete optimization problem with bounded submodularity ratio, leading to polynomial-time algorithms with provable optimality bounds. We show that our approach is applicable to different types of uncertainties, including additive and multiplicative uncertainties in the system matrices as well as uncertain time delays. We demonstrate our approach in a numerical case study on the IEEE 39-bus test power system.

## I. INTRODUCTION

Networked systems are often controlled by providing control signals to a subset of input nodes, which then influence the rest of the network via local interactions. Since controlling a large number of input nodes is often impractical, selecting a minimum-size set of input nodes to satisfy desired system properties has been an area of recent research interest. Algorithms have been proposed for selecting input nodes based on criteria including controllability [1], [2], stability [3], [4], and robustness to noise [5], as well as application-specific methodologies for power systems [6].

Existing input selection studies consider either linear systems with known, fixed parameters, systems with a discrete set of possible parameter values, or systems with arbitrary real-valued parameters [1], [2], [3], [4]. Typically, however, systems will have structured uncertainties of different types, such as small additive or multiplicative perturbations of matrix entries. Inputs that are selected based on a given nominal system model have been shown to fail to guarantee stability in the presence of such perturbations [7]. Moreover, current works assume that system states are not affected by delays, which may further degrade stability.

In this paper, we study the problem of selecting a minimum-size input set to guarantee existence of a stabilizing feedback control for a linear time-invariant system with delays and uncertainties. Our approach is based on two insights. First, we show that the existence of such a controller is guaranteed if and only if a subset of modes, which is determined by the $H_\infty$ norm of the uncertainty, lies in the controllability subspace induced by the set of input nodes. Second, we show that the problem of selecting a minimum-cardinality set of input nodes that satisfies this condition is equivalent to a column-subset selection problem [8]. We make the following specific contributions:

- We derive sufficient conditions for a set of input nodes to guarantee existence of a stabilizing feedback controller in an uncertain dynamical system. Our conditions are based on the set of system modes that are within the controllability subspace resulting from the input nodes, and are expressed as functions of the magnitude of the uncertainty. We extend our results to stability of time-delayed systems by deriving conditions for Padé approximations to time delays.
- We formulate the problem of selecting a minimum-size set of input nodes to ensure robustness to a given type and magnitude of uncertainty, and prove that this problem is equivalent to minimizing the Euclidean distance between the identified system modes and the controllability subspace. We prove that this Euclidean distance is a monotone decreasing function of the input set with bounded submodularity ratio.
- We propose polynomial-time approximation algorithms for selecting a minimum-cardinality input set guaranteeing system robustness. We characterize the optimality gap/bound of our approach via the submodularity ratio, which we derive as a function of the system parameters. We analyze different types of uncertainties under our framework, including additive and multiplicative uncertainties, as well as uncertain time delays in outputs.
- We evaluate our input selection approach through a numerical study on the IEEE 39-bus test power system [9]. The goal is to select a minimal set of generators to participate in the wide-area damping control to ensure the system stability. Our results show that the proposed approach better ensures the system reachability in presence of uncertainties compared to existing techniques that do not take uncertainties into account. By comparison to a geometric index-based selection scheme, our approach requires fewer control inputs to achieve the system robustness condition.

The paper is organized as follows. Section II reviews the related work. Section III introduces the system model. Section IV presents the proposed approach to input selec-

Z. Liu, Y. Long, P. Lee, L. Bushnell, D. Kirschen, and R. Poovendran are with the Department of Electrical Engineering, University of Washington, Seattle, WA 98195 USA. {zhipliu, longyao, leep3, lb2, kirschen, rp3}@uw.edu

A. Clark is with the Department of Electrical and Computer Engineering, Worcester Polytechnic Institute, Worcester, MA 01609 USA. aclark@wpi.edu

This work was supported by NSF grant CNS-1544173.

tion for robust control. Section V contains analysis of the proposed approach for different uncertainties. Section VI presents numerical results. Section VII concludes the paper.

## II. RELATED WORK

Input selection for control systems has received extensive research attention in recent years. In [10], the problem of identifying a minimum-size set of input nodes to control a directed network was studied. The problem of computing the sparsest diagonal input matrix for a given LTI system that guarantees the system controllability, known as the minimal controllability problem, is addressed in [1], and greedy algorithms are proposed to approximate the optimal solution with provable guarantees on the sparsity. The exact solution to the minimal controllability problem is explored in [2] which shows that this problem can be reduced to the minimum set covering problem.

Input selection with different inputs and constraints has been considered in [11], [12], [13], [14], [5]. In [11], the authors study the problem of minimizing total control effort for a given state transfer while ensuring controllability. In [12] and [13], the problem constraint is to ensure structural controllability, defined as existence of a controllable numerical realization of the linear system matrices $(A, B)$ with the same structure (i.e., zero/nonzero pattern). The submodularity properties of the minimal input selection in dynamical networks are explored in [14] and [5], in which efficient selection algorithms are proposed that achieve solutions with provable optimality bounds. All these works focus on controllability of the system, while our approach ensures the system stability, defined as the ability of a system to reach to an equilibrium.

In many control problems, controllability is a sufficient but not necessary condition to ensure system stability. While controllability implies that all system modes are controllable, stabilizability only requires that unstable modes be controllable [4]. The problem of selecting minimal inputs to achieve stabilizability is addressed in [3], based on a metric related to the span of the controllability matrix. A submodular framework for input selection to ensure stabilizability was presented in [4] along with applications to small signal stability of power systems. These works assume known system parameters and hence may lose the stability guarantee in a system with plant uncertainties or time delays.

Uncertainties and time delays are common problems in any LTI systems. The problem of selecting inputs to guarantee the design of a robust control in the presence of uncertainties and delays, however, has not been studied in the existing literature. Initial discussions on input selection for robust control are proposed in [7] which is based on a random searching strategy that does not guarantee the feasiblity or optimality of the solutions. In this paper, we propose an optimization approach using bounded submodularity ratio to minimal input selection problem that guarantees the existence of a robust control and achieves system reachability from a given state to the origin.

The hardness of approximating minimal reachability problems is presented in [15], in which it is proved that there is no polynomial-time algorithm to achieve a set of inputs within a constant factor of true optimal selections. In this paper, we show that the bounded submodularity ratio alone is sufficient to achieve certain approximations, even without submodularity or supermodularity. The optimality bound, however, can be arbitrarily small for certain system parameters and hence is not contradicting with [15].

## III. LINEAR SYSTEM MODEL WITH UNCERTAINTIES

In this section, we present a generalized linear system model with uncertainties and a transformation to the $M - \Delta$ loop representation [16]. We also review controllability metrics, submodularity ratio, and linear regressions, and introduce notations that will be used throughout the paper.

We consider a linear time-invariant dynamical system of $n$ states with a perturbation $\delta(t)$, described by

$$\dot{x}(t) = Ax(t) + Bu(t) + \delta(t) \quad (1)$$
$$y(t) = Cx(t) \quad (2)$$

where $x(t) \in \mathbb{R}^n$ is the state vector; $y(t) \in \mathbb{R}^m$ is a vector of $m$ outputs; $u(t) \in \mathbb{R}^p$ is the input vector whose $i$th entry is denoted by $u_i$ that represents the $i$th input signal into the system. The system matrix $A \in \mathbb{R}^{n \times n}$ has eigenvalues $\{\lambda_i\}$ and corresponding eigenvectors $\{v_i\}$. The input matrix $B \in \mathbb{R}^{n \times p}$ has columns $\{b_1, \ldots, b_p\}$ where each column $b_i$ corresponds to the influence of a unit control input $u_i$ to the system. The output matrix $C \in \mathbb{R}^{m \times n}$ has rows $\{c_1^T, \ldots, c_m^T\}$ where each row $c_i^T$ represents the connection between an output $y_i$ and the states $x$.

We consider a feedback control, given as follows, based on the state estimation

$$u(t) = -K\hat{x}(t), \quad (3)$$

where $K$ is the feedback matrix to be designed and $\hat{x}(t)$ is an estimate of the system state $x(t)$ with dynamics

$$\dot{\hat{x}}(t) = A\hat{x}(t) + Bu(t) + L(y(t) - C\hat{x}(t)).$$

In state estimation, $L$ is a weighting matrix on the output error $y(t) - C\hat{x}(t)$ that is used to correct the estimation dynamics. Defining the estimation error $e(t) = x(t) - \hat{x}(t)$ induces the following dynamics

$$\dot{e}(t) = (A - LC)e(t) + \delta(t). \quad (4)$$

By involving the feedback control (3) and state estimation dynamics (4), the closed-loop system dynamics of Eq. (1) becomes

$$\begin{bmatrix} \dot{x} \\ \dot{e} \end{bmatrix} = \begin{bmatrix} A - BK & BK \\ 0 & A - LC \end{bmatrix} \begin{bmatrix} x \\ e \end{bmatrix} + \begin{bmatrix} I \\ I \end{bmatrix} \delta. \quad (5)$$

The uncertainty vector $\delta(t)$ can be decomposed based on the specific type of uncertainty. For example, when a perturbation occurs in matrix $A$ and results in an uncertain system $\dot{x}(t) = (A + \Delta)x(t) + Bu(t)$, the uncertainty term is represented by $\delta(t) = \Delta x(t)$.

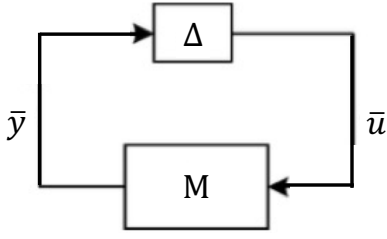

Fig. 1: M-$\Delta$ loop for stability analysis

Given uncertainties in terms of the matrix $\Delta$ and the state $x$, the closed-loop system (5) can be transformed into an $M-\Delta$ loop representation for stability analysis [16], described by the following equations

$$\dot{\overline{x}} = A_{cl}\overline{x} + B_{cl}\overline{u} \tag{6}$$
$$\overline{y} = C_{cl}\overline{x} \tag{7}$$
$$\overline{u} = \Delta \overline{y} \tag{8}$$

The new variables relate to those in Eq. (5) by the following equations

$$\overline{x} = \begin{bmatrix} x \\ e \end{bmatrix}, \quad A_{cl} = \begin{bmatrix} A - BK & BK \\ 0 & A - LC \end{bmatrix},$$

$$B_{cl}\Delta C_{cl}\overline{x} = \begin{bmatrix} I \\ I \end{bmatrix}\delta.$$

The parameters $B_{cl}$ and $C_{cl}$ are constant matrices whose values depend on the decomposition of $\delta$. For example, when $\delta(t) = \Delta x(t)$, we have $B_{cl} = \begin{bmatrix} I & I \end{bmatrix}^T$ and $C_{cl} = \begin{bmatrix} I & 0 \end{bmatrix}$. Denote $b$ and $c$ as the smallest numbers satisfying $B_{cl}B_{cl}^T \preceq bI$ and $C_{cl}^T C_{cl} \preceq cI$. The transformation to $M-\Delta$ system for different types of uncertainties will be studied in Section V.

Fig. 1 shows a block diagram of the $M-\Delta$ loop, where the block $\Delta$ is the uncertainty matrix while the block $M$ contains all certain dynamics of the system. The transfer function of the certain block $M$ is given by

$$M = \left[\begin{array}{c|c} A_{cl} & B_{cl} \\ \hline C_{cl} & 0 \end{array}\right].$$

In what follows, we review some definitions related to the controllability, submodularity ratio, and linear regressions.

The controllability matrix for a linear system $\dot{x} = Ax + Bu$ is defined by $\mathcal{C} = [B \ AB \ A^2B \ \ldots \ A^{n-1}B]$.

*Definition 1 (Submodular ratio [8]):* Let $f : 2^\Omega \to \mathbb{R}$ be a non-negative set function. The submodularity ratio of $f$ with respect to a set $U \subseteq \Omega$ and a parameter $k \geq 1$ is

$$\gamma_{U,k} = \min_{\substack{L,S: L \cap S = \emptyset \\ L \subseteq U, |S| \leq k}} \frac{\sum_{x \in S}(f(L \cup \{x\}) - f(L))}{f(L \cup S) - f(L)}$$

The submodularity ratio measures to what extent $f$ has submodular properties. For a nondecreasing function $f$, we have $\gamma_{U,k} \in [0,1]$ and $f$ is submodular if and only if $\gamma_{U,k} \geq 1$.

*Definition 2 (Coefficient of determination [8]):* For a linear regression problem $\min_x \|Hx - z\|_2^2$, assuming the columns of $H \in \mathbb{R}^{n \times m}$ and $z \in \mathbb{R}^n$ are normalized to have norm 1, the coefficient of determination (or $R^2$ statistic) is

$$R^2 = n \ \overline{b}^T \overline{C}^{-1} \overline{b}$$

where $\overline{C} = H^T H/n$ and $\overline{b} = H^T z/n$.

## IV. PROPOSED INPUT SELECTION FRAMEWORK

This section formulates the problem of selecting a minimal subset of inputs that guarantees the existence of a robust feedback control in the presence of system uncertainties. We first give the exact formulation. In order to achieve computational feasibility, we then derive a sufficient condition. We prove that selecting a minimum-size set of input nodes to satisfy the sufficient condition is equivalent to a discrete optimization problem with bounded submodularity ratio. By exploiting the bounded submodularity ratio, we present a polynomial-time greedy algorithm with a provable optimality guarantee on the size of selection.

### A. Problem Formulation

Given a set of input nodes, $\Omega = \{1, \ldots, p\}$, selecting a subset $S \in \Omega$ to exert control signals changes the system dynamics (1) to

$$\dot{x}(t) = Ax(t) + B_s u_s(t) + \delta(t) \tag{9}$$

where the input vector $u_s$ has entries $\{u_i : i \in S\}$ and the matrix $B_s$ consists of columns $\{b_i : i \in S\}$. The input selection therefore impacts the system dynamics by determining the columns of the matrix $B$.

The goal is to select the minimal set of inputs that guarantees the existence of a feedback control, $u_s = -K\hat{x}$ where $\hat{x}$ is the state estimate, robust to the uncertainty $\delta$.

In the equivalent $M-\Delta$ loop system (6)-(8) of dynamics (9), the selection $S$ affects the matrix $B_s$ and hence changes the matrix $A_{cl}$ to

$$A_{cl} = \begin{bmatrix} A - B_s K & B_s K \\ 0 & A - LC \end{bmatrix}. \tag{10}$$

We assume the uncertainty matrix $\Delta$ in an $M-\Delta$ system has bounded singular values, i.e., $\|\Delta\|_\infty \leq \sigma$. The following lemma gives an equivalent condition on the selection $S$ for existence of a robust control that stabilizes the system in the presence of uncertainties $\Delta$.

*Lemma 1 ([16], Small Gain Theorem):* The interconnected system in Fig. 1 is well posed and internally stable for all $\Delta$ with $\|\Delta\|_\infty \leq \sigma$ if and only if $\|M\|_\infty < 1/\sigma$.

The following lemma shows an equivalent condition to the robustness requirement given by Lemma 1.

*Lemma 2 ([17], Lemma 7.4):* Suppose $M = \left[\begin{array}{c|c} A_{cl} & B_{cl} \\ \hline C_{cl} & 0 \end{array}\right]$. Then $\|M\|_\infty < 1/\sigma$ if and only if there exists a positive definite matrix $X \succ 0$ such that

$$\begin{bmatrix} \sqrt{\sigma}C_{cl}^T \\ 0 \end{bmatrix}\begin{bmatrix} \sqrt{\sigma}C_{cl} & 0 \end{bmatrix} + \begin{bmatrix} (A_{cl}^T X + X A_{cl}) & \sqrt{\sigma}X B_{cl} \\ \sqrt{\sigma}B_{cl}^T X & -I \end{bmatrix} \prec 0. \tag{11}$$

By Schur complement theorem, the inequality (11) is equivalent to

$$\sigma C_{cl}^T C_{cl} + A_{cl}^T X + X A_{cl} + \sigma X B_{cl} B_{cl}^T X \prec 0. \quad (12)$$

By Lemma 1 and 2, the problem of selecting a minimal set of inputs to ensure the existence of a feedback control robust to uncertainties $\|\Delta\|_\infty \leq \sigma$, can be formulated as

$$\min_S \ |S| \quad (13)$$
$$\text{s.t.} \ \sigma C_{cl}^T C_{cl} + A_{cl}^T X + X A_{cl} + \sigma X B_{cl} B_{cl}^T X \prec 0 \quad (14)$$
$$\text{for some} \ X \succ 0$$

Note that (14) contains the constructed matrix $A_{cl}$ defined by Eq. (10) and hence is a function of the selection $S$.

Without any additional structure on the constraint (14), solving the problem (13)-(14) is difficult as it requires exhaustive search over all $S$. In what follows, we present a sufficient condition to satisfy the constraint of (14) such that the problem (13)-(14) is relaxed to a problem that does not require information of matrix $X$.

*Lemma 3:* Suppose $B_{cl} B_{cl}^T \preceq bI$ and $C_{cl}^T C_{cl} \preceq cI$. Suppose that there exists $\epsilon > 0$ such that $X \prec \epsilon I$ when $X$ is defined by

$$X = \int_0^\infty e^{A_{cl}^T t}(\sigma c I + \sigma b \epsilon^2 I) e^{A_{cl} t} \, dt. \quad (15)$$

Then X is a solution to (14).
*Proof:* Suppose $B_{cl} B_{cl}^T \preceq bI, \ C_{cl}^T C_{cl} \preceq cI$. Then

$$\sigma C_{cl}^T C_{cl} + A_{cl}^T X + X A_{cl} + \sigma X B_{cl} B_{cl}^T X$$
$$\preceq \sigma c I + A_{cl}^T X + X A_{cl} + \sigma b X^2.$$

Thus, any positive definite matrix $X$ that can solve

$$\sigma c I + A_{cl}^T X + X A_{cl} + \sigma b X^2 \prec 0 \quad (16)$$

is a solution to inequality (14).

The constructed $X$ in Eq. (15) can be verified to be a solution of (16) by substituting it into inequality (16). When $X \prec \epsilon I$, it can be shown that $X^2 \prec \epsilon^2 I$ (see appendix for proof) and hence we have

$$\sigma c I + A_{cl}^T X + X A_{cl} + \sigma b X^2$$
$$= \sigma c I + \sigma b X^2 + \int_0^\infty A_{cl}^T e^{A_{cl}^T t}(\sigma c I + \sigma b \epsilon^2 I) e^{A_{cl} t} \, dt$$
$$+ \int_0^\infty e^{A_{cl}^T t}(\sigma c I + \sigma b \epsilon^2 I) e^{A_{cl} t} A_{cl} \, dt$$
$$= \sigma c I + \sigma b X^2 + \int_0^\infty \frac{d}{dt}\left(e^{A_{cl}^T t}(\sigma c I + \sigma b \epsilon^2 I) e^{A_{cl} t}\right) dt$$
$$= \sigma c I + \sigma b X^2 + e^{A_{cl}^T t}(\sigma c I + \sigma b \epsilon^2 I) e^{A_{cl} t}\Big|_{t=0}^\infty$$
$$= \sigma c I + \sigma b X^2 - (\sigma c I + \sigma b \epsilon^2 I)$$
$$= \sigma b (X^2 - \epsilon^2 I) \prec 0$$

which implies that the constructed $X$ in (15) is a solution to inequality (16) and hence is also a solution to (14). ∎

By Lemma 3, in order to meet the constraint (14) of the minimal input selection problem, it suffices to require the matrix $A_{cl}$ in Eq. (10) to satisfy that for some $\epsilon > 0$,

$$\int_0^\infty e^{A_{cl}^T t}(\sigma c I + \sigma b \epsilon^2 I) e^{A_{cl} t} \, dt \prec \epsilon I. \quad (17)$$

The following lemma simplifies the inequality (17), which results in a sufficient condition to constraint (14).

*Lemma 4:* Let $\lambda_0$ denote the eigenvalue of $A_{cl}$ with the largest real part. Define $\alpha = \lambda_1(A_{cl} + A_{cl}^T)/(2\text{Re}(\lambda_0))$, where the notation $\lambda_1(\cdot)$ refers to the largest eigenvalue. Assume $A_{cl}$ is asymptotically stable, i.e., $\text{Re}(\lambda_0) < 0$. Then the inequality (17) holds if

$$\text{Re}(\lambda_0) < -\sigma\sqrt{bc}/\alpha. \quad (18)$$

*Proof:* Given the fact that a symmetric matrix $(\epsilon I - X)$ is positive definite if and only if all its eigenvalues are positive, the inequality (17) holds if the largest eigenvalue of $X$ is smaller than $\epsilon$. The constructed $X$ in (17) is a solution of the Lyapunov equation

$$A_{cl}^T X + X A_{cl} + Q = 0, \quad \text{where} \ Q = (\sigma c + \sigma b \epsilon^2) I.$$

The eigenvalues of solution $X$ are bounded by $\lambda_1(X) \leq -\lambda_1(Q)/\lambda_1(A_{cl} + A_{cl}^T)$ when $A_{cl} + A_{cl}^T \prec 0$ [18]. Replacing $\lambda_1(A_{cl} + A_{cl}^T)$ by $2\alpha\text{Re}(\lambda_0)$, we have

$$\lambda_1(X) \leq \frac{\sigma c + \sigma b \epsilon^2}{-2\alpha \text{Re}(\lambda_0)}.$$

To guarantee $\lambda_1(X) < \epsilon$, it suffices to require

$$\text{Re}(\lambda_0) < -\frac{\sigma}{2\alpha}\left(\frac{c}{\epsilon} + b\epsilon\right). \quad (19)$$

After finding the minimum value of the right-hand-side function in (19) respect to $\epsilon$, we have the inequality (19) holds for some $\epsilon > 0$ if $\text{Re}(\lambda_0) < -\sigma\sqrt{bc}/\alpha$. ∎

By Lemma 3 and 4, we can relax the constraint of problem (13)-(14) to a sufficient condition (18) on the eigenvalues of matrix $A_{cl}$, which still guarantees the design of a robust feedback control $K$ that can stabilize the system with uncertainty $\Delta$. Hence, the problem of selecting minimal inputs to ensure the existence of a robust control in presence of bounded uncertainties $\|\Delta\|_\infty \leq \sigma$ can be formulated as

$$\min_S \ |S| \quad (20)$$
$$\text{s.t.} \ \text{Re}(\lambda_0) < -\sigma\sqrt{bc}/\alpha \quad (21)$$

where $\lambda_0$ is the eigenvalue of $A_{cl}$ with largest real part and $B_{cl} B_{cl}^T \preceq bI, \ C_{cl}^T C_{cl} \preceq cI$.

By the fact that $2\text{Re}(\lambda_0) \leq \lambda_1(A_{cl} + A_{cl}^T) < 0$ [19], we have that $\alpha$ is in the range $(0, 1]$. We take the following approach to approximating the solution to (21). We first set $\alpha = 1$, yielding the constraint

$$\text{Re}(\lambda_0) < -\sigma\sqrt{bc}, \quad (22)$$

and select a set $S$ using the procedure enumerating in the following sections. After synthesizing a controller and computing the matrices $A_{cl}$, $B_{cl}$, and $C_{cl}$, if the condition of Lemma 2 does not hold, then we reduce the value of $\alpha$ and use the procedure of the following sections, but with the constraint $Re(\lambda_0) < -\sigma\sqrt{bc}/\alpha$. We continue this procedure until an input set $S$ satisfying Lemma 2 is selected. In practice, during the simulations described in Section VI, the condition (22) was sufficient to ensure robust stability during all simulation trials.

### B. Input Selection Constraints

The matrix $A_{cl}$ in (10) possesses block upper triangular structure, which maps the eigenvalue constraint (22) to

$$\text{Re}(\lambda) < -\sigma\sqrt{bc} \quad (23)$$

for all eigenvalues $\lambda$ of $A - B_s K$ and $A - LC$.

Assume the state estimation matrix $L$ is chosen such that $A - LC$ has all eigenvalues satisfying the condition (23). Then it suffices to require $A - B_s K$ to satisfy (23), in order to achieve a stable system in the presence of uncertainties satisfying $\|\Delta\|_\infty \leq \sigma$. In other words, if the selection $S$ can guarantee the existence of a feedback control $K$ such that $A - B_s K$ has all eigenvalues smaller than $-\sigma\sqrt{bc}$, then the controlled system is stable and robust to any bounded uncertainties $\|\Delta\|_\infty \leq \sigma$.

*Lemma 5:* Let

$$\dot{x}(t) = Ax(t) + B_s u(t) \quad (24)$$
$$y(t) = Cx(t) \quad (25)$$

be a control system. If all of the eigenvectors of $A$ with eigenvalues $\lambda$ satisfying $\text{Re}(\lambda) \geq \lambda_1$ lie in the span of the controllability matrix of the system, then there exists a feedback control $K$ such that all eigenvalues of the closed-loop system $A - B_s K$ satisfy $\text{Re}(\lambda) < \lambda_1$.

*Proof:* The proof can be found in appendix, which is a straightforward generalization of Chen [20] and is included for completeness. ∎

Lemma 5 maps the input selection constraint (23) to a condition that requires all eigenvectors of "undesired modes" ($\text{Re}(\lambda) \geq -\sigma\sqrt{bc}$) to be contained in the span of the controllability matrix of the system $(A, B_s)$, i.e.,

$$\{v_i : \text{Re}(\lambda_i) \geq -\sigma\sqrt{bc}\} \in \text{span}(\mathcal{C}(S)), \quad (26)$$

where $\mathcal{C}(S) = [B_s \ AB_s \ A^2 B_s \ \ldots \ A^{n-1} B_s]$. Those "undesired modes" with eigenvalues $\text{Re}(\lambda) \geq -\sigma\sqrt{bc}$ can be viewed as unstable modes in the presence of bounded uncertainties $\|\Delta\|_\infty \leq \sigma$.

It can be shown that for the linear system (24) with selected inputs $S$, the span of the controllability matrix is equal to the span of the controllability Gramian [21], i.e.,

$$\text{span}(\mathcal{C}(S)) = \text{span}(W(S)), \quad (27)$$

where the controllability Gramian of the system $(A, B_s)$ is defined by

$$W(S) = \int_{t_0}^{t_1} e^{A(t-t_0)} B_s B_s^T e^{A^T(t-t_0)} dt$$

for some $t_1 > t_0$.

Motivated by the input selection constraint (26) and (27), we define the metric $F(S)$ by

$$F(S) \triangleq \sum_{i: \ Re(\lambda_i) \geq -\sigma\sqrt{bc}} dist^2(v_i, \text{span}(W(S))), \quad (28)$$

where $dist(\cdot)$ denotes the Euclidean distance. The metric $F(S)$ can be interpreted as the distance between the eigenvectors of "undesired modes" where $\text{Re}(\lambda) \geq -\sigma\sqrt{bc}$ and the span of the controllability Gramian. Intuitively, this metric is a measure of how close the unstable modes in the presence of uncertainty $\Delta$ are to being controllable.

By Lemma 5 and the discussion in (26) and (27), the constraint $F(S) = 0$ is sufficient to ensure the existence of a robust feedback control that can stabilize the system in presence of bounded uncertainties $\|\Delta\|_\infty \leq \sigma$.

### C. Submodularity Ratio of Input Selection Constraint

In this subsection, we prove that the function $F(S)$ defined in (28) has bounded submodularity ratio. Before giving the bound, we first define some notations as follows.

Let $P \in \mathbb{R}^{np \times np}$ be a nonsingular matrix that normalizes the controllability matrix $\mathcal{C}$ to $\mathcal{C}' = \mathcal{C}P$ where each column of $\mathcal{C}'$ has norm 1. For simplicity, we use $\mathcal{C}_s$ to denote the function $\mathcal{C}(S)$ in the following. Denote $\mathcal{C}'_s = \mathcal{C}_s P$.

Define $\bar{C} \in \mathbb{R}^{np \times np} = \mathcal{C}'^T \mathcal{C}'/n$ and $\bar{C}_s = \mathcal{C}'^T_s \mathcal{C}'_s/n$. We note that $\bar{C}_s$ is a submatrix[1] of $\bar{C}$ with rows and columns selected by set $S$. For any matrix $\bar{C}_s$, we denote its smallest eigenvalue as $\lambda_{\min}(\bar{C}_s)$ and let $\lambda_{\min}(\bar{C}, k) = \min_{S: |S|=k} \lambda_{\min}(\bar{C}_s)$.

*Theorem 1:* Let $\gamma_{U,k}$ be the submodularity ratio of the set function $F(S)$ defined in (28). Then for any set $U \subseteq \Omega$ and $k \geq 1$, the submodularity ratio $\gamma_{U,k}$ is bounded by

$$\gamma_{U,k} \geq \lambda_{\min}(\bar{C}, k + |U|) \geq \lambda_{\min}(\bar{C}).$$

Before giving the proof, we first present some preliminary results as follows.

Given any vector $v$ with $\|v\|_2 = 1$, consider the function

$$f_v(S) = dist^2(v, \text{span}(\mathcal{C}(S))) = \min_x \|\mathcal{C}(S)x - v\|_2^2.$$

*Lemma 6:* For any nonsingular matrix $P$, the following equality holds: $\min_x \|\mathcal{C}_s x - v\|_2^2 = \min_x \|\mathcal{C}_s Px - v\|_2^2$.

*Proof:* Let $x'^* = \arg\min_x \|\mathcal{C}_s Px - v\|_2^2$. Then $x'^* = P^{-1}(\mathcal{C}_s^T \mathcal{C}_s)^{-1} P^{-T} P^T \mathcal{C}_s^T v = P^{-1}(\mathcal{C}_s^T \mathcal{C}_s)^{-1} \mathcal{C}_s^T v = P^{-1} x^*$ where $x^* = \arg\min_x \|\mathcal{C}_s x - v\|_2^2$. Thus,

$$\min_x \|\mathcal{C}_s Px - v\|_2^2 = \|\mathcal{C}_s Px'^* - v\|_2^2 = \|\mathcal{C}_s PP^{-1} x^* - v\|_2^2$$
$$= \|\mathcal{C}_s x^* - v\|_2^2 = \min_x \|\mathcal{C}_s x - v\|_2^2,$$

completing the proof. ∎

By Lemma 6, the function $f(S)$ is equivalent to

$$f_v(S) = \min_x \|\mathcal{C}'_s x - v\|_2^2. \quad (29)$$

---
[1]$\bar{C}_s$ consists of $n \times |S|$ rows and columns selected from $\bar{C}$.

Denote $x_s^*$ as the optimal vector that solves $\min_x \|\mathcal{C}_s' x - v\|_2^2$. The value of $x_s^*$ is given by $x_s^* = (\mathcal{C}_s'^T \mathcal{C}_s')^{-1} \mathcal{C}_s'^T v$ [22]. Known that $v - \mathcal{C}_s' x_s^*$ and $\mathcal{C}_s' x_s^*$ are orthogonal, we have

$$f_v(S) = \|\mathcal{C}_s' x_s^* - v\|_2^2 = \|v\|_2^2 - \|\mathcal{C}_s' x_s^*\|_2^2$$
$$= 1 - \|\mathcal{C}_s' x_s^*\|_2^2.$$

Define a new set function

$$g_v(S) = \|\mathcal{C}_s' x_s^*\|_2^2 = \|P_s v\|_2^2, \tag{30}$$

where $P_s = \mathcal{C}_s'(\mathcal{C}_s'^T \mathcal{C}_s')^{-1} \mathcal{C}_s'^T$ is the projection matrix for orthogonal projection onto the span of columns of $\mathcal{C}_s'$, and hence

$$g_v(S) = 1 - f_v(S). \tag{31}$$

*Lemma 7 ([23]):* The function $g_v(S)$ is the $R^2$ statistic (coefficient of determination) for the linear regression problem (29).

*Proof:* By definition,

$$g_v(S) = \|\mathcal{C}_s'(\mathcal{C}_s'^T \mathcal{C}_s')^{-1} \mathcal{C}_s'^T v\|_2^2$$
$$= (\mathcal{C}_s'^T v)^T (\mathcal{C}_s'^T \mathcal{C}_s')^{-1} (\mathcal{C}_s'^T v) = n\, \bar{v}_s^T \bar{C}_s^{-1} \bar{v}_s,$$

where $\bar{C}_s = \mathcal{C}_s'^T \mathcal{C}_s'/n$ and $\bar{v}_s = \mathcal{C}_s'^T v/n$, satisfying the definition of $R^2$ statistic for the linear regression (29) [8]. ∎

*Lemma 8 ([8]):* Let $\gamma_{U,k}'$ be the submodularity ratio of the $R^2$ statistic (30). Then

$$\gamma_{U,k}' \geq \lambda_{\min}(\bar{C}, k + |U|) \geq \lambda_{\min}(\bar{C}). \tag{32}$$

By Lemma 7 and 8, the function $g_v(S)$ has submodularity ratio $\gamma_{U,k}'$ bounded by (32). Using this result, now we are ready to give the proof of Theorem 1.

*Proof of Theorem 1:* We first show that $f_v(S)$ and $g_v(S)$ have the same submodularity ratio $\gamma_{U,k}'$. From (31), we have $g_v(S) = 1 - f_v(S)$. By definition,

$$\gamma_{U,k}' = \min_{\substack{L,S: L\cap S = \emptyset \\ L \subseteq U, |S| \leq k}} \frac{\sum_{x \in S} (g_v(L \cup \{x\}) - g_v(L))}{g_v(L \cup S) - g_v(L)}$$

$$= \min_{\substack{L,S: L\cap S = \emptyset \\ L \subseteq U, |S| \leq k}} \frac{\sum_{x \in S} ((1 - f_v(L \cup \{x\})) - (1 - f_v(L)))}{(1 - f_v(L \cup S)) - (1 - f_v(L))}$$

$$= \min_{\substack{L,S: L\cap S = \emptyset \\ L \subseteq U, |S| \leq k}} \frac{\sum_{x \in S} (f_v(L \cup \{x\}) - f_v(L))}{f_v(L \cup S) - f_v(L)}$$

Then we show that $F(S)$ has submodularity ratio $\gamma_{U,k}$ bounded by $\gamma_{U,k}'$. By Eq. (27) and Lemma 6, we have $F(S) = \sum_i f_{v_i}(S)$ where $v_i$ is the eigenvector of the $i$th "undesired mode" with eigenvalue $(\text{Re}(\lambda) \geq -\sigma\sqrt{bc})$. Then

$$\gamma_{U,k} = \min_{\substack{L,S: L\cap S = \emptyset \\ L \subseteq U, |S| \leq k}} \frac{\sum_{x \in S} (F(L \cup \{x\}) - F(L))}{F(L \cup S) - F(L)}$$

$$= \min_{\substack{L,S: L\cap S = \emptyset \\ L \subseteq U, |S| \leq k}} \frac{\sum_i \sum_{x \in S} (f_{v_i}(L \cup \{x\}) - f_{v_i}(L))}{\sum_i (f_{v_i}(L \cup S) - f_{v_i}(L))}$$

$$\geq \min_{\substack{L,S: L\cap S = \emptyset \\ L \subseteq U, |S| \leq k}} \left( \min_i \frac{\sum_{x \in S} (f_{v_i}(L \cup \{x\}) - f_{v_i}(L))}{(f_{v_i}(L \cup S) - f_{v_i}(L))} \right)$$

$$\geq \gamma_{U,k}'$$

**Algorithm 1** Algorithm for selecting minimal set of inputs for robust control.

1: **procedure** MINInput($\Omega$)
2:    **Input:** Set of input nodes $\Omega$
3:    $S \leftarrow \emptyset$
4:    **while** $F(S) > 0$ and $S \neq \Omega$ **do**
5:      **for** $v \in \Omega \setminus S$ **do**
6:        $\delta_v \leftarrow F(S) - F(S \cup \{v\})$
7:      **end for**
8:      $v^* \leftarrow \arg\max_v \delta_v$
9:      $S \leftarrow S \cup \{v^*\}$
10:    **end while**
11:    **return** $S$
12: **end procedure**

By Lemma 8, we have $\gamma_{U,k} \geq \gamma_{U,k}' \geq \lambda_{\min}(\bar{C}, k + |U|) \geq \lambda_{\min}(\bar{C})$, completing the proof. ∎

### D. Minimal Input Selection Algorithm

Using the metric established in the previous section, the problem of selecting the minimal set of inputs to ensure the existence of a robust control in presence of uncertainties can be approached by the following formulation:

$$\min_S |S|, \quad s.t.\ F(S) = 0. \tag{33}$$

The bounded submodularity ratio of $F(S)$ implies that there exists a polynomial-time greedy algorithm to approximately solve it with provable optimality bound when the submodularity ratio is not zero [8], [23]. In what follows, we present such a greedy algorithm described by Algorithm 1.

The algorithm proceeds as follows. The set of selected inputs, $S$, is initialized to be empty. At each iteration, the element $v \in \Omega \setminus S$ that maximizes $F(S) - F(S \cup \{v\})$ is selected and added to $S$. The algorithm terminates when $F(S)$ reaches zero or $S = \Omega$. Note that if $F(\Omega) > 0$, there is no guarantee on the system stability for any set of inputs.

The optimality bound of Algorithm 1 is defined by the following proposition.

*Proposition 1:* Let $S^*$ denote the true optimal solution to problem (33) and let $S$ denote the solution returned by Algorithm 1. Then we have

$$\frac{|S| - 1}{|S^*|} \leq \frac{1}{\lambda_{\min}(\bar{C}, 2|S|)} \log \frac{F(\emptyset)}{F(S_{t-1})},$$

where $S_{t-1}$ denotes the selected input set $S$ at the second-to-last iteration of Algorithm 1.

*Proof:* Let $\{s_1, \ldots, s_k\}$ denote the elements selected by the greedy algorithm in the first $k$ iterations. By the definition of submodularity ratio and Theorem 1, we have

$$\frac{\sum_{x \in S^*} (F(\emptyset) - F(\{x\}))}{F(\emptyset) - F(S^*)} \geq \gamma_{S, |S|} \geq \gamma_0$$

where $\gamma_0 = \lambda_{\min}(\bar{C}, 2|S|)$, which implies that

$$|S^*| (F(\emptyset) - F(\{s_1\})) \geq \gamma_0 (F(\emptyset) - F(S^*)),$$

or equivalently,
$$F(\{s_1\}) - F(S^*) \leq \left(1 - \frac{\gamma_0}{|S^*|}\right)(F(\emptyset) - F(S^*)).$$

Similarly, we have that
$$\frac{\sum_{x \in S^*}(F(s_1) - F(\{x, s_1\}))}{F(\{s_1\}) - F(S^* \cup \{s_1\})} \geq \gamma_0$$

implying that
$$F(\{s_1\}) - F(\{s_1, s_2\}) \geq \frac{\gamma_0}{|S^*|}(F(\{s_1\}) - F(S^* \cup \{s_1\}))$$
$$\geq \frac{\gamma_0}{|S^*|}(F(\{s_1\}) - F(S^*)),$$

which is further equivalent to
$$F(\{s_1, s_2\}) - F(S^*) \leq \left(1 - \frac{\gamma_0}{|S^*|}\right)(F(\{s_1\}) - F(S^*))$$
$$\leq \left(1 - \frac{\gamma_0}{|S^*|}\right)^2 (F(\{\emptyset\}) - F(S^*)).$$

By mathematical induction, we have that
$$F(\{s_1, \ldots, s_k\}) - F(S^*) \leq \left(1 - \frac{\gamma_0}{|S^*|}\right)^k (F(\{\emptyset\}) - F(S^*)).$$

Note that $F(S^*) = 0$ and hence
$$F(\{s_1, \ldots, s_k\}) \leq \left(1 - \frac{\gamma_0}{|S^*|}\right)^k F(\{\emptyset\}).$$

Use the fact that $1 - \frac{1}{x} \leq \log(x), \forall x \geq 1$. Then when
$$k = \frac{|S^*|}{\gamma_0} \log \frac{F(\emptyset)}{F(S_{t-1})},$$

we have
$$F(\{s_1, \ldots, s_k\}) \leq F(S_{t-1}),$$

implying that one additional element is sufficient. Hence
$$\frac{|S| - 1}{|S^*|} \leq \frac{1}{\gamma_0} \log \frac{F(\emptyset)}{F(S_{t-1})}.$$
∎

For a set $\Omega = \{1, \ldots, p\}$, Algorithm 1 terminates in at most $p$ iterations. Each iteration solves at most $n$ least-squares problems where each problem needs computation of $O(n^3)$, which gives an overall complexity of $O(pn^4)$.

We note that the optimality bound in Proposition 1 can be arbitrarily small for certain system matrices. This is consistent with the claim that minimal reachability is hard to approximate as presented in [15].

## V. ANALYSIS FOR DIFFERENT UNCERTAINTIES

This section presents a detailed discussion on the $M - \Delta$ loop transformation of the linear system (1)-(2) for different uncertainties $\delta$. We also derive the specific values of $b$ and $c$ in the metric $F(S)$ for each uncertainty case.

### A. Additive Uncertainty in Matrix A

When an uncertain matrix $\Delta$ is additive to matrix $A$, the linear system has dynamics $\dot{x} = (A + \Delta)x + Bu$, $y = Cx$. When mapping this system to (1)-(2), the uncertainty vector in (1) has decomposition $\delta = \Delta x$. In the $M - \Delta$ transformation (6)-(8), we have $B_{cl} = \begin{bmatrix} I \\ I \end{bmatrix}$ and $C_{cl} = \begin{bmatrix} I & 0 \end{bmatrix}$.

In this case, we get $b = 2$ and $c = 1$ using the fact that
$$B_{cl}B_{cl}^T = \begin{bmatrix} I & I \\ I & I \end{bmatrix} \preceq 2\begin{bmatrix} I & 0 \\ 0 & I \end{bmatrix}, \quad C_{cl}^T C_{cl} = \begin{bmatrix} I & 0 \\ 0 & 0 \end{bmatrix} \preceq \begin{bmatrix} I & 0 \\ 0 & I \end{bmatrix}$$

Thus, for additive uncertainties in matrix $A$, the metric (28) of the input selection problem for robust control becomes
$$F(S) = \sum_{i:\ Re(\lambda_i) \geq -\sigma\sqrt{2}} dist^2(v_i, \text{span}(W(S))). \quad (34)$$

### B. Multiplicative Uncertainty in Matrix A

A linear system with multiplicative uncertainties in matrix $A$ has dynamics $\dot{x} = (I + \Delta)Ax + Bu$, $y = Cx$.

Following analysis of the previous case, the uncertainty vector can be written as $\delta = \Delta Ax$ and hence in the $M - \Delta$ transformation, $B_{cl} = \begin{bmatrix} I \\ I \end{bmatrix}$ and $C_{cl} = \begin{bmatrix} A & 0 \end{bmatrix}$.

For this case, we also have $b = 2$ by the same fact, $B_{cl}B_{cl}^T \preceq 2I$, as in previous case. For the parameter $c$, let $\rho$ denote the largest singular value of $A$. Then we have
$$C_{cl}^T C_{cl} = \begin{bmatrix} A^T A & 0 \\ 0 & 0 \end{bmatrix} \preceq \rho^2 \begin{bmatrix} I & 0 \\ 0 & I \end{bmatrix},$$

which implies $c = \rho^2$. Hence the metric $F(S)$ for input selection in presence of multiplicative uncertainties in matrix $A$ is given by
$$F(S) = \sum_{i:\ Re(\lambda_i) \geq -\sigma\rho\sqrt{2}} dist^2(v_i, \text{span}(W(S))).$$

### C. Uncertain Time Delay in Output y

For a time delay $\tau$ in output $y$, we denote the delayed output by $y_d(t) = y(t - \tau)$. In frequency domain, the delay can be modeled by a first-order Pade approximation [24],
$$y_d = e^{-s\tau}y \approx \frac{-0.5s\tau + 1}{0.5s\tau + 1}y. \quad (35)$$

A state-space representation of (35) is given by
$$\dot{x}_d(t) = \frac{2}{\tau}(-x_d(t) + 2y(t)) \quad (36)$$
$$y_d(t) = x_d(t) - y(t) \quad (37)$$

When each output $y_i$ has a different time delay $\tau_i$, define a diagonal matrix of delays, $\Gamma = diag(\tau_1, \ldots, \tau_m)$, and $n$ equations of (36)-(37) are added to the linear system $\dot{x} = Ax + Bu$, $y = Cx$, which gives the overall system dynamics,
$$\begin{bmatrix} \dot{x} \\ \dot{x}_d \end{bmatrix} = \begin{bmatrix} A & 0 \\ 4\Gamma^{-1}C & -2\Gamma^{-1} \end{bmatrix} \begin{bmatrix} x \\ x_d \end{bmatrix} + \begin{bmatrix} B \\ 0 \end{bmatrix} u \quad (38)$$
$$y_d = \begin{bmatrix} -C & I \end{bmatrix} \begin{bmatrix} x \\ x_d \end{bmatrix} \quad (39)$$

Suppose the uncertain delay variables $\Gamma^{-1} = \Gamma_0^{-1} + \Delta$, where $\Gamma_0^{-1}$ is the known certain delays and $\Delta$ is the uncertain delays with the assumption that $\|\Delta\|_\infty \leq \sigma$. By defining a new state vector $\tilde{x} = [x, x_d]^T$, the system (38)-(39) can be re-written in the form of (1)-(2):

$$\dot{\tilde{x}} = \tilde{A}\tilde{x} + \tilde{B}u + \tilde{\delta} \qquad (40)$$
$$y_d = \tilde{C}\tilde{x} \qquad (41)$$

where $\tilde{A} = \begin{bmatrix} A & 0 \\ 4\Gamma_0^{-1}C & -2\Gamma_0^{-1} \end{bmatrix}$, $\tilde{B} = \begin{bmatrix} B \\ 0 \end{bmatrix}$, $\tilde{C} = \begin{bmatrix} -C & I \end{bmatrix}$ and $\tilde{\delta} = \begin{bmatrix} 0 & 0 \\ 4\Delta C & -2\Delta \end{bmatrix} \tilde{x}$.

When putting the system (40)-(41) into the $M - \Delta$ representation, we have

$$A_{cl} = \begin{bmatrix} \tilde{A} - \tilde{B}K & \tilde{B}K \\ 0 & \tilde{A} - L\tilde{C} \end{bmatrix},$$
$$B_{cl} = \begin{bmatrix} 0 & 4I & 0 & 4I \end{bmatrix}^T, \quad C_{cl} = \begin{bmatrix} I & -\tfrac{1}{2}I & 0 & 0 \end{bmatrix}.$$

Studying the largest eigenvalues of $B_{cl}$ and $C_{cl}$ gives $B_{cl}B_{cl}^T \preceq 32I$ and $C_{cl}^T C_{cl} \preceq 5/4$. Therefore we have $b = 32$ and $c = 5/4$, and the input selection metric is given by

$$F(S) = \sum_{i: \ Re(\lambda_i) \geq -2\sigma\sqrt{10}} dist^2(v_i, \text{span}(W(S))).$$

## VI. NUMERICAL STUDIES

This section presents simulation results for the proposed input selection algorithm (Algorithm 1) in wide-area damping control of power systems [25]. By linearizing around a steady state, the generator angle swings have the state-space form of dynamics (1)-(2). The goal of the wide-area damping control is to damp the inter-area oscillations between generators and stabilize all generators' rotor speeds. The control input signals are adjusting voltages at each generator to impact the power output, while observations (control outputs) consist of generator power outputs and transmission line power flows. In this control system, uncertainties may arise due to load changing, unexpected transmission line tripping and time delay.

The problem is to select a set of generators to be involved in the wide-area damping control while guaranteeing that a robust control exists to damp generator rotor swings in presence of uncertainties. We test our minimal input selection approach on the IEEE 39-bus test system [9], and compare it with the case when the metric $f(S)$ does not take uncertainties into account and a geometric indices based selection method [26].

The system topology and initial parameters including generator settings and transmission line impedances are specified in [9]. We obtain the initial operating point with one unstable mode by scaling all power system stabilizers' gain parameters to 16.2% and increasing overall load to 105% of initial level. Linearizing the system dynamics around current operating point gives a linear system with matrices $A, B, C$ in the state-space form (1)-(2). For any input selection $S$, we have the corresponding linear system $(A, B_s, C)$ and we use the same linear-quadratic regulator (LQR) design to compute the feedback control $u = -Kx$, which minimizes the following quadratic cost function,

$$J(u) = \int_0^\infty (100 \ x^T x + u^T u) \ dt.$$

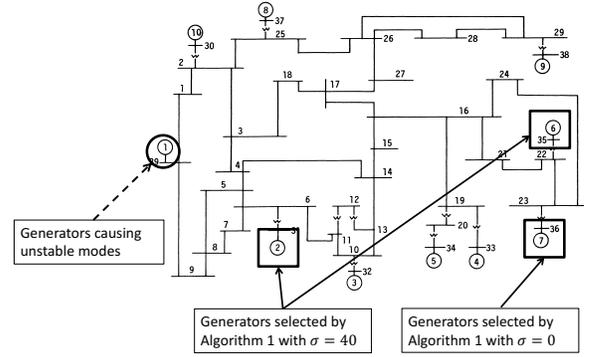

Fig. 2: IEEE 10-generator 39-bus power system topology. Circled generators have initial rotor angles differing from the rest, which causes unstable modes. Generators selected by our proposed algorithm with thresholds $\sigma = 0$ and $\sigma = 40$ are indicated by squares.

We consider scenarios when a random matrix $\Delta$ is added to matrix $A$. The metric (34) is used in the proposed Algorithm 1. Assuming the largest singular value of $\Delta$ is less than 40, our proposed algorithm selects generators $\{2, 6\}$ to participate in the wide-area control.

To evaluate how the assumption $\|\Delta\|_\infty \leq \sigma$ on uncertainties impacts our decision, we compare two selections using the same metric (34) with different thresholds $\sigma = 0$ and $\sigma = 40$ respectively. When $\sigma = 0$, the selected inputs ensure stability of the nominal system but do not guarantee robustness to uncertainties, while $\sigma = 40$ results in a selection for robust control assuming $\|\Delta\|_\infty \leq 40$. The reason for choosing the selection with $\sigma = 0$ as a baseline is that to the best of our knowledge, there is no other analytical work on input selection to ensure system robustness in the existing literature. In this test case, the Algorithm 1 selects a single generator $\{7\}$ when $\sigma = 0$ and generators $\{2, 6\}$ when $\sigma = 40$. Fig. 2 shows a comparison of selections for $\sigma = 0$ and $\sigma = 40$ on the power system topology.

Fig. 3 shows a comparison of robustness between the two systems with inputs selected by $\sigma = 0$ and $\sigma = 40$. Each data point is the percentage of cases that the controlled system is stable in 1000 trials of random matrix $\Delta$ with the same largest singular value. We observe that when the uncertainty matrix $\Delta$ has largest singular value within the assumption $\|\Delta\|_\infty \leq 40$, the feedback control based on the selection of $\sigma = 40$ ensured stability in all cases.

When the uncertainty matrix $\Delta$ has singular values over the bound, for example $\|\Delta\|_\infty = 76.59$, both selections fail to guarantee stability. The algorithm with robust consideration $\sigma = 40$, however, achieves a system with more stable cases (91.1%) compared to the selection based on $\sigma = 0$ (76.0%) in 1000 random trials.

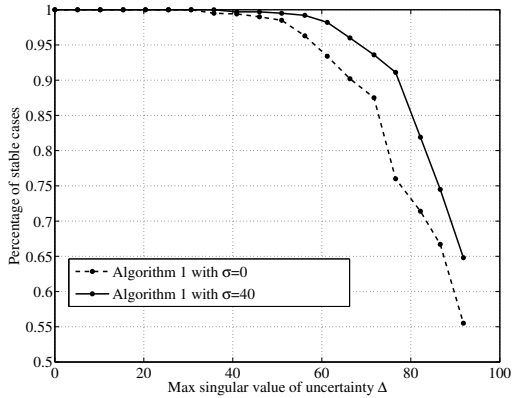

Fig. 3: Comparison of percentages of cases that the controlled system is stable in the presence of uncertainty $\Delta$ using the proposed greedy selection with $\sigma = 0$ and the selection with $\sigma = 40$. Each data point has 1000 random trials in $\Delta$. The proposed greedy algorithm with robustness consideration ($\sigma = 40$) is stable in a higher percentage of cases compared to the greedy selection without considering uncertainties ($\sigma = 0$).

A comparison between our proposed Algorithm 1 and the geometric indices based selection [26] is shown in Table I. The geometric selection is implemented as follows. Assuming the uncertaintiy is bounded by $\|\Delta\|_\infty \leq \sigma$, we first compute the geometric indices of controllability $m_{ci}(k)$ for each generator $i$ and each undesired mode $k$ with $\mathrm{Re}(\lambda) \geq -\sigma\sqrt{2}$,

$$m_{ci}(k) = \frac{(b_i^T \psi_k)}{(\|\psi_k\|_2 \|b_i\|_2)},$$

where $b_i$ is the $i$th column of matrix $B$ and $\psi_k$ is the left eigenvector of matrix $A$ corresponding to the eigenvalue $\lambda_k$. Then the algorithm selects the generator with the highest controllability index in each iteration until the metric (34) reaches zero.

We observe that in order to guarantee the system stability in presence of bounded uncertainties, our proposed input selection requires fewer generators compared to the geometric indices based selection.

TABLE I: Comparing generators selected by proposed greedy algorithm and geometric algorithm.

| $\sigma$ | 0 | 0.01 | 0.1 | 0.3 | 0.6 | 0.7 | 1 |
|---|---|---|---|---|---|---|---|
| Proposed selection | {7} | {8} | {5, 9} | {2, 5} | {3, 4} | {3, 6} | {2, 6} |
| Geometric selection | {9} | {10} | {9, 10} | {9, 10} | {4, 9, 10} | {4, 7, 9} | {4, 8, 9} |

## VII. CONCLUSIONS

In this paper, we studied the problem of selecting the minimal set of inputs for robust control in order to ensure system stability in presence of uncertainties. We formulated a sufficient condition for the input selection to guarantee the existence of a stabilizing control, requiring that a subset of system modes lie in the span of the controllability Gramian. We exploited the bounded submodularity ratio of the minimal input selection problem and presented a computationally efficient algorithm that has provable optimality guarantees. We also derived specific metrics used in the proposed algorithm for additive system uncertainties, multiplicative uncertainties and output time delays respectively. We evaluated our approach in a case study on the IEEE 39-bus power system with a comparison to existing methods. Our future work will include output selection for robust control and joint design of input selection and robust control in order to minimize the total control effort.

## ACKNOWLEDGMENTS

The authors would like to thank A. Jadbabaie, A. Olshevsky, G. Pappas, and V. Tzoumas for pointing out an error in an earlier version of Theorem 1 and its proof [27, Lemma 6].

## APPENDIX

*Lemma 9:* For any symmetric positive definite matrix $X$, if $X \prec \epsilon I$, then $X^2 \prec \epsilon^2 I$.

*Proof:* For any $X \succ 0$, we have $\epsilon^2 I - X^2 = (\epsilon I - X)(\epsilon I + X)$, where $\epsilon I - X$ is positive definite by assumption that $X \prec \epsilon I$ and $\epsilon I + X \succ X \succ 0$ is also positive definite.

The product of two positive definite matrices $A$ and $B$ is positive definite if $A$ and $B$ commute, i.e., $AB = BA$. Let $A = \epsilon I - X$ and $B = \epsilon I + X$. Hence we have $AB = \epsilon^2 I - X^2$ positive definite, which implies $X^2 \prec \epsilon^2 I$. ∎

*Proof of Lemma 5:* Let $\{q_1, \ldots, q_{n_1}\}$ be a maximal set of linearly independent columns of the controllability matrix $\mathcal{C}$, and define a matrix $P$ such that $P^{-1} = (q_1 \cdots q_{n_1} \cdots q_n)$, where $q_{n_1+1}, \ldots, q_n$ can be arbitrarily chosen so long as $P^{-1}$ is nonsingular. Under the transformation $z = Px$, the system can be written as $\dot{z} = \overline{A}z + \overline{B}u$ or more specifically

$$\begin{bmatrix} \dot{z}_1(t) \\ \dot{z}_2(t) \end{bmatrix} = \begin{bmatrix} A_{11} & A_{12} \\ 0 & A_{22} \end{bmatrix} \begin{bmatrix} z_1 \\ z_2 \end{bmatrix} + \begin{bmatrix} B_1 \\ 0 \end{bmatrix} u(t),$$

where the $n_1$-dimensional subsystem $\dot{z}_1 = A_{11}z_1 + B_1 u$ is controllable (Chen [20, Theorem 6.6]).

Suppose first that $\lambda$ is an eigenvalue where the corresponding eigenvector $v$ lies in the span of the controllability matrix, i.e., $v \in span(\mathcal{C}) = span([B_s \ AB_s \ A^2 B_s \ \ldots \ A^{n-1} B_s])$. We show that $\lambda$ is an eigenvalue of $A_{11}$, and hence all eigenvalues with $\text{Re}(\lambda) \geq \lambda_1$ are eigenvalues of $A_{11}$.

Indeed, letting $w = Pv$, we have that $\overline{A}w = PAP^{-1}w = PAP^{-1}Pv = PAv = \lambda Pv = \lambda w$, establishing that $w$ is an eigenvector of $\overline{A}$ corresponding to $\lambda$. Furthermore, we have that $v = P^{-1}w$, and hence $w$ is the representation of $v$ with respect to the columns of $P^{-1}$. Since $v$ is in the span of $\mathcal{C}$, it also lies in the span of $\{q_1, \ldots, q_{n_1}\}$, and thus $w_i = 0$ for $i > n_1$. We then have that $w = (w_1 \ 0)^T$, where $w_1$ is an eigenvector of $A_{11}$ corresponding to eigenvalue $\lambda$.

A similar argument to the above establishes that all eigenvalues corresponding to "undesired modes" ($\text{Re}(\lambda) \geq \lambda_1$) are eigenvalues of $A_{11}$, and that the system $(A_{11}, B_1)$ is controllable. Hence we can design a feedback controller $K$ such that the system $A_{11} - B_1 K$ has eigenvalues in the desired region ($\text{Re}(\lambda) < \lambda_1$), while the remaining eigenvalues are unchanged ($\text{Re}(\lambda) < \lambda_1$). Applying the transformation $P$ and using the fact that the eigenvalues of the system are invariant under a similarity transform yields the desired result. ∎